\documentclass[11pt]{article}

\usepackage{fullpage}

\usepackage{amssymb}
\usepackage{amsthm}
\usepackage[english]{babel}
\usepackage{url}
\usepackage{algorithm}
\usepackage{algorithmic}
\usepackage{tabularx}
\usepackage{paralist}
\usepackage{mathtools}
\usepackage{tcolorbox}
\usepackage{hyperref}
\usepackage[capitalise]{cleveref}
\usepackage{caption}
\usepackage{subcaption}

\usepackage{tcolorbox}

\usepackage{nicefrac} 

\usepackage{colortbl}
\definecolor{bgcolor}{rgb}{0.8, 1, 1}
\definecolor{bgcolor2}{rgb}{0.8, 1, 0.8}
\definecolor{niceblue}{rgb}{0.0, 0.19, 0.56}

\hypersetup{
    colorlinks=true,
    linkcolor=blue,
    citecolor=niceblue,
    urlcolor=blue
}

\newtheorem{theorem}{Theorem}
\newtheorem{proposition}{Proposition}
\newtheorem{corollary}{Corollary}

\newtheorem{assumption}{Assumption}

\newtheorem{remark}{Remark}
\newtheorem{lemma}{Lemma}

\renewcommand{\|}{\parallel}

\newcommand\ev[1]{\left \langle #1 \right \rangle}

\newcommand{\norm}[1]{\left\lVert#1\right\rVert}
\newcommand{\sqn}[1]{{\left\lVert#1\right\rVert}^2}

\usepackage{pifont}
\definecolor{PineGreen}{RGB}{0,110,51}
\definecolor{BrickRed}{RGB}{143,20,2}

\usepackage{tikz-cd} 

\usepackage{xspace}

\makeatletter

\usepackage[T1]{fontenc}

\newcommand{\mytag}[1]{%
  \refstepcounter{equation}%
  \edef\@currentlabel{\theequation}%
  {({\@currentlabel})}%
  \@bsphack
  \begingroup
    \@onelevel@sanitize\@currentlabelname
    \edef\@currentlabelname{%
      \expandafter\strip@period\@currentlabelname\relax.\relax\@@@%
    }%
    \protected@write\@auxout{}{%
      \string\newlabel{#1}{%
        {\@currentlabel}%
        {\thepage}%
        {\@currentlabelname}%
        {\@currentHref}{}%
      }%
    }%
  \endgroup
  \@esphack
}
\makeatother

\title{\bf Incorporating Preconditioning into Accelerated Approaches: Theoretical Guarantees and Practical Improvement}
\author{\begin{tabular}{c}
     {\bf Stepan Trifonov}\\
     {\small MCAS}
\end{tabular}\quad \begin{tabular}{c}
     {\bf Leonid Levin}\\
     {\small MCAS} 
\end{tabular} \quad  \begin{tabular}{c}
     {\bf Savelii Chezhegov} \\
     {\small MCAS}\thanks{Corresponding author: \texttt{savachezhegov2017@gmail.com}.} 
\end{tabular} \\
\begin{tabular}{c}
     {\bf Aleksandr Beznosikov}\\
     {\small MCAS}
\end{tabular}}
\begin{document}
\maketitle
\begin{abstract}

Machine learning and deep learning are widely researched fields that provide solutions to many modern problems. Due to the complexity of new problems related to the size of datasets, efficient approaches are obligatory. In optimization theory, the Heavy Ball and Nesterov methods use \textit{momentum} in their updates of model weights. On the other hand, the minimization problems considered may be poorly conditioned, which affects the applicability and effectiveness of the aforementioned techniques. One solution to this issue is \textit{preconditioning}, which has already been investigated in approaches such as \textsc{AdaGrad}, \textsc{RMSProp}, \textsc{Adam} and others. Despite this, momentum acceleration and preconditioning have not been fully explored together. Therefore, we propose the Preconditioned Heavy Ball (\textsc{PHB}) and Preconditioned Nesterov method (\textsc{PN}) with theoretical guarantees of convergence under \textit{unified} assumption on the scaling matrix. Furthermore, we provide numerical experiments that demonstrate superior performance compared to the unscaled techniques in terms of iteration and oracle complexities.

\end{abstract}

\section{Introduction}
A classical optimization problem can be stated as
\begin{align}
\label{eq: problem}
    \min_{x \in \mathbb{R}^d} f(x),
\end{align}
where the function $f$ can have different meanings depending on the setting of the machine learning \cite{shalev2014understanding} and deep learning \cite{goodfellow2016deep} problems, such as loss function, risk function, etc. The most intuitive, yet simple, approach for solving the problem \eqref{eq: problem} is Gradient Descent \cite{cauchy1847methode}. However, there are other first-order oracle approaches allow for both theoretical and practical improvements. These methods include the Heavy-Ball method \cite{polyak1964} and Accelerated Gradient Descent \cite{nesterov1983}, which serve as the foundation for a wide range of approaches aimed at addressing narrowly focused problems.  

Although the momentum technique provides acceleration, it does not address ill-conditioned problems. Therefore, techniques related to the scaling of directions in the iterative process are relevant here. First, it is worth mentioning quasi-Newton methods \cite{fletcher1963rapidly,goldfarb1970family,shanno1970conditioning,liu1989limited} based on the Hessian approximation. The essence of this approach is quite simple: instead of using the exact inverse Hessian in Newton's method, we use its approximation for efficiency, as calculating the inverse matrix of second derivatives is currently regarded as an expensive oracle. 

Subsequently, approaches that are categorized as \textit{adaptive} began to emerge due to their use of update rules that pull information from previous points during the iterative process. Such approaches include \textsc{AdaGrad} \cite{duchi2011}, \textsc{RMSProp} \cite{tieleman2012lecture}, \textsc{Adam} \cite{kingma2014}, \textsc{AMSGrad} \cite{reddi2019convergence}, \textsc{NAdam} \cite{dozat2016incorporating}, \textsc{OASIS} \cite{jahani2021doubly}, \textsc{AdaHessian} \cite{yao2021adahessian} and many others. Their distinctive feature from Newton's method and its approximations is that the preconditioning matrix has a diagonal form (as a consequence, it becomes much cheaper to reverse it). Along with their empirical superiority over other existing techniques, adaptive methods are frequently chosen to train machine learning models. 

At present, considerable research is devoted to the convergence analysis of adaptive methods, as well as to the integration of preconditioning into techniques such as variance reduction \cite{sadiev2024stochastic}, extragradient \cite{beznosikov2022scaled} and distributed approaches \cite{reddi2020adaptive,chezhegov2024local} in a generalized form. Despite this, scaling has not been considered under general assumptions on the preconditioning matrix for accelerated approaches.\\\\
\noindent \textbf{Our contribution.} Therefore, our contribution can be formulated as follows:
\begin{itemize}
    \item We present two algorithms that are scaled versions of the Heavy-Ball and Nesterov methods.
    \item We propose theoretical guarantees of presented methods under unified assumption on the precondioning matrix and validate numerical experiments to show the outperformance compared to unscaled techniques.
\end{itemize}
\section{Related Works}
Currently, methods such as the Heavy-Ball method and Accelerated Gradient Descent have been investigated in detail under various assumptions concerning the target function and, in some cases, on the first-order stochastic oracle \cite{polyak1964,nesterov1983,ghadimi2012optimal,flammarion2015averaging}. In this paper, we rely on the analysis of the Heavy-Ball method from \cite{danilova2022convergence} and the classical analysis of Nesterov three-point method \cite{nesterov2013introductory}.


If we discuss preconditioning techniques, we should mention methods that use only gradient, such as \textsc{Adagrad} \cite{duchi2011}, which originally allowed to solve the problem of nonsmooth optimization. Subsequent developments extended this approach through exponential moving averages (e.g., \textsc{RMSProp} \cite{tieleman2012lecture}) and momentum techniques (e.g., \textsc{Adam} \cite{kingma2014}), now widely adopted as default optimizers in deep learning. Moreover, several variants of \textsc{Adam} have since been proposed, including Nesterov momentum -- \textsc{NAdam} \cite{dozat2016incorporating}, modified second moment estimation -- \textsc{AMSGrad} \cite{reddi2019convergence}, and
decoupled weight decay -- \textsc{AdamW} \cite{loshchilov2017decoupled}.

There are also approaches where a Hessian rather than the gradient itself is used to adapt the scaling matrix. For example, such methods include \textsc{OASIS} \cite{jahani2021doubly}, \textsc{SOPHIA} \cite{liu2023sophia} and \textsc{AdaHessian} \cite{yao2021adahessian} -- methods with Hutchinson approximation which aim to leverage curvature information while retaining scalability.

As already mentioned, preconditioning has already been considered under the general assumption as a modification of existing approaches. For example, variance reduced methods such as \textsc{Scaled SARAH} and \textsc{Scaled L-SVRG} \cite{sadiev2024stochastic} was derived. Moreover, for the saddle point problem, a scaled version of the extragradient method \cite{beznosikov2022scaled} was considered. In the distributed optimization framework, the scaled \textsc{Local SGD} method \cite{chezhegov2024local} was already analyzed.

\section{Preliminaries, requirements and notations}
Throughout this paper, we adopt the following notations. The scalar product in $L_2$-space is denoted as $\ev{x, y}$. We denote $\norm{x} = \norm{x}_2 = \sqrt{\ev{x, y}}$ as the norm in $L_2$-space.
The induced scalar product is denoted as $\ev{x, y}_{A} = \ev{x, Ay}$. The norm induced by a positive definite matrix $A$ we denote as $\norm{x}_A = \sqrt{\ev{x, Ax}}$. The Hadamard product of two same-dimensional elements $A$ and $B$ is denoted as $A \odot B$.
\subsection{Objective function}
\noindent As for the assumptions on the target function, we suppose that the minimized function $f$ from \eqref{eq: problem} is $\mu$-strongly convex and $L$-smooth.
\begin{assumption}[$\mu$-strong convexity] 
\label{asm: strong}
The function $f$ is $\mu$-strongly convex, i.e $\forall x, y \in \mathbb{R}^d$
\begin{align*}
    f(y) \geq f(x) + \ev{\nabla f(x), y - x} + \frac{\mu}{2}\norm{y - x}^2.
\end{align*}   
\end{assumption}
\begin{assumption}[$L$-smoothness] 
\label{asm: smooth}
The function $f$ is $L$-smooth, i.e $\forall x, y \in \mathbb{R}^d$
\begin{align*}
    \norm{\nabla f(y) - \nabla f(x)} \leq L\norm{y - x}.
\end{align*}
\end{assumption}
\noindent Also $L$-smoothness can be written in an equivalent formulation which is more useful for theoretical proofs.
\begin{proposition}
\label{prop: equiv-L}
Assume that $f(x)$ satisfies \cref{asm: smooth}. Therefore, it holds that $\forall x \in \mathbb{R}^d$
\begin{align*}
    \norm{\nabla f(x)}^2 \leq 2L(f(x) - f(x^*)).
\end{align*}
\end{proposition}
\subsection{Preconditioning}
\label{sec: prec}
The update of the scaling matrix may differ from one approach to another. Nevertheless, this procedure can be described as sequence of diagonal matrices $\{D_{k}\}_{k=-1}^{\infty}$, where each $D_k$ can be evaluated by the following rule:
\begin{align}
\label{eq: scale-update}
    D_k = \text{Update}(D_{k-1}, H_k),
\end{align}
where $H_k$ is some information received at iteration $k$ in the diagonal form. For instance, \textsc{AdaGrad} \cite{duchi2011} update can be written as 
\begin{align*}
    D_k = \sqrt{D_{k-1}^2 + \text{diag}(\nabla f(x_k) \odot \nabla f(x_k))}.
\end{align*}
Moreover, some methods use so-called \textit{exponential smoothing}. Prominent examples of such approaches are \textsc{RMSProp} \cite{tieleman2012lecture} and \textsc{Adam} \cite{kingma2014}, where the update \eqref{eq: scale-update} can be written as
\begin{align}
\label{eq: quadr-update-adam}
    D_k = \sqrt{\beta_2 D_{k-1}^2 + (1 - \beta_2)\text{diag}(\nabla f(x_k) \odot \nabla f(x_k))},
\end{align}
where the parameter $\beta_2$ lies in 
$(0,1)$. In theoretical studies, the standard choice for the smoothing parameter $\beta_2$ is $\beta_2(k) = 1 - \frac{1}{k}$ or $\beta_2 = 1 - \frac{1}{K}$, while in practice a commonly used value is $0.999$.

\noindent Sometimes, instead of exponential squared smoothing, exponential smoothing for linear summands can be used, such as in the \textsc{OASIS} \cite{jahani2021doubly} method:
\begin{align}
\label{eq: linear-update-oasis}
    D_k = \beta_2 D_{k-1} + (1 - \beta_2)\text{diag}(v_k \odot \nabla^2 f(x_k)v_k),
\end{align}
where $v_k$ is the vector of random variables sampled from the Rademacher distribution.\\
In general form, \eqref{eq: quadr-update-adam}  and \eqref{eq: linear-update-oasis} can be represented as
\begin{align}
\label{eq: quadr-up}
    D_k = \sqrt{\beta_2 D_{k-1}^2 + (1 - \beta_2)H_k^2},
\end{align}
or
\begin{align}
\label{eq: linear-up}
    D_k = \beta_2 D_{k-1} + (1 - \beta_2)H_k.
\end{align}
We devote our analysis exactly to these rules, \eqref{eq: quadr-up}-\eqref{eq: linear-up}, for updating the matrix $D_k$. A natural assumption that is made on the sequence $\{D_k\}_{k = -1}^\infty$ is the boundedness of the scaling matrices.
\begin{assumption}
\label{asm: scaling-matr}
The sequence $\{D_k\}_{k=-1}^\infty$ is bounded, i.e. there exist constants $\Gamma \geq e > 0$ such that for all $k$
\begin{align*}
    eI \preceq D_k \preceq \Gamma I.
\end{align*} 
\end{assumption}
\noindent This assumption is indeed valid in terms of $\Gamma$ -- usually this value is strongly correlated with the dimensionality of the problem and the smoothness constant (\cite{chezhegov2024local}, Table 1). However, the assumption of boundedness from below leaves questions. In this case, the following trick can be applied:
\begin{align}
\label{eq: trick}
    [\hat{D}_k]_{ii} := \max\{e, [D_k]_{ii}\}.
\end{align}
Let us also formulate two auxiliary facts about preconditioning matrices.
\begin{proposition}
\label{prop: scale-bound}
Assume that the matrix $A$ satisfies
\begin{align*}
    eI \preceq A \preceq \Gamma I
\end{align*}
for some constants $\Gamma \geq e > 0$. Therefore,
\begin{align*}
    e\norm{x}^2 &\leq \norm{x}^2_{A} \leq \Gamma\norm{x}^2;\\
    \frac{1}{\Gamma}\norm{x}^2 &\leq \norm{x}^2_{A^{-1}} \leq \frac{1}{e}\norm{x}^2.
\end{align*}
\end{proposition}
To prove \cref{prop: scale-bound}, it is enough to apply the definition of the induced norm.
\begin{proposition} [\cite{beznosikov2022scaled} -- Lemma 1, \cite{chezhegov2024local} -- Corollary 1]
\label{prop: ultra}
Assume that the sequence of matrices $\{D_{k}\}_{k=-1}^{\infty}$ satisfies \cref{asm: scaling-matr}. If the update rule \eqref{eq: scale-update} has the form \eqref{eq: quadr-up}/\eqref{eq: linear-up} with \eqref{eq: trick}, then the following inequality holds:    
\begin{align*}
    \norm{x}_{D_{k+1}}^2 \leq (1 + (1 - \beta_2)C)\norm{x}_{D_{k}}^2,
\end{align*}
where the constant $C$ depends on the preconditioner update rule. To be more precise, 
\begin{align*}
    C = \begin{cases}
        \frac{\Gamma^2}{2e^2} &\qquad \text{for \eqref{eq: quadr-up}};\\
        \frac{2\Gamma}{e} &\qquad \text{for \eqref{eq: linear-up}}.
    \end{cases}
\end{align*}
\end{proposition}

\section{Algorithms and Theoretical Results}
In this section, we present the design of \textsc{Preconditioned Heavy-Ball Method} and \textsc{Preconditioned Nesterov Method} with its theoretical guarantees of convergence. 
\begin{algorithm}
\caption{\textsc{Preconditioned Heavy-Ball Method}}
\label{alg:scaling_algorithm}
\begin{algorithmic}[1]
\REQUIRE{initial point $x_0$, parameter of momentum $\beta_1$, initial momentum $V_{-1} = 0$, initial scaling matrix $\hat{D}_{-1} \succeq eI$}
\FOR{$k = 0, 1, 2, \ldots, T-1$}
    \STATE $V_{k}\;=\;\beta_1 V_{k-1}\;+\;\hat{D}_{k}^{-1}\nabla f(x_{k})$
    \STATE $x_{k+1}\;=\;x_{k}\;-\;\gamma V_{k}$
    \STATE $[\hat{D}_{k + 1}]_{ii} = \max\{e, [\text{Update}(D_{k}, H_{k+1})]_{ii}\}$
\ENDFOR
\end{algorithmic}
\end{algorithm}

\begin{algorithm}
\caption{\textsc{Preconditioned Nesterov Method}}
\label{alg:scaling_algorithm_nesterov}
\begin{algorithmic}[1]
\REQUIRE{initial point $x^0 = x^0_f = x^0_g$, learning rate $\gamma_k$, initial scaling matrix $\hat{D}_{-1} \succeq eI$, momentums $\{\xi_k\}_{k=0} \geq 0, \{\theta_k\}_{k=0} \geq 0$}
\FOR{$k = 0, 1, 2, \ldots, T-1$}
    \STATE $x^{k + 1}_f = x^k_g - \gamma_k \hat{D}_k^{-1} \nabla f(x^k_g)$
    \STATE $x^{k + 1} = \xi_k (x^{k + 1}_f - x_{f}) + x_f^{k}$
    \STATE $x_g^{k + 1} = \theta_{k + 1} x_{f}^{k + 1} + (1 - \theta_{k + 1}) x^{k + 1}$
    \STATE $[\hat{D}_{k + 1}]_{ii} = \max\{e, [\text{Update}(D_{k}, H_{k+1})]_{ii}\}$
\ENDFOR
\end{algorithmic}
\end{algorithm}
\noindent The algorithmic design is fairly standard. Unlike classical Gradient Descent, we employ either the heavy-ball momentum or Nesterov momentum as an acceleration mechanism. At the same time, the core modification of existing accelerated methods lies in the scaling matrix $\hat{D}_k$, which is updated according to rules \eqref{eq: quadr-up}/\eqref{eq: linear-up} and \eqref{eq: trick}. As previously mentioned, many adaptive methods naturally fit into such preconditioning schemes, making our proposed algorithms, in a sense, unified frameworks.

\noindent We now formulate the theorems that provide convergence guarantees for Algorithms \ref{alg:scaling_algorithm} and \ref{alg:scaling_algorithm_nesterov}.
\begin{tcolorbox}[colback=gray!25, colframe=gray!10, sharp corners]
\begin{theorem}
\label{thm:1}
Suppose that Assumptions \ref{asm: strong}, \ref{asm: smooth} and \ref{asm: scaling-matr}  hold. Then, after $K$ iterations of \cref{alg:scaling_algorithm} with $\gamma = \frac{(1 - \beta_1)^2e}{12 L}$, we have
\begin{align*}
     f\left(\frac{1}{W_{K-1}}\sum\limits_{k=0}^{K-1} w_kx_k\right) - f(x^*) &\leq 4\exp\left(-\frac{(1 - \beta_1)\mu e K}{48 L \Gamma}\right) L\norm{{x}_{0} - x^{*}}_{\hat{D}_{-1}}^{2},
\end{align*}
where $W_{K-1} = \sum_{k=0}^{K-1} w_k$, $w_k = \left(1 - \frac{\mu F}{4 \Gamma}\right)^{-(k+1)}$.
\end{theorem}
\end{tcolorbox}
\begin{remark}
    It is worth noting that the technique of point reweighting is standard in the analysis of scaling-based methods. For further details, see Appendix.
\end{remark}
As a consequence, the upper bound on the required number of iterations $K$ for reaching $\varepsilon$-accuracy can be formulated as follows.
\begin{tcolorbox}[colback=gray!25, colframe=gray!10, sharp corners]
\begin{corollary}
     Under the conditions of \cref{thm:1}, the required number of iterations $K$ of \cref{alg:scaling_algorithm} for reaching $\varepsilon$-accuracy, i.e. $f(x_\text{out}) -f(x^*) \leq \varepsilon$, can be upper bounded as
    \begin{align*}
        K = \mathcal{O}\left(\frac{L\Gamma}{\mu e (1 - \beta_1)}\log\left(\frac{L\sqn{x_0 - x^*}}{\varepsilon}\right)\right).
    \end{align*}
\end{corollary}
\end{tcolorbox}
As for the \cref{alg:scaling_algorithm_nesterov}, the theorem of convergence can be formulated as follows.
\begin{tcolorbox}[colback=gray!25, colframe=gray!10, sharp corners]
\begin{theorem}
\label{thm:2}
Suppose that Assumptions \ref{asm: strong}, \ref{asm: smooth} and \ref{asm: scaling-matr}  hold. Then, after $K$ iterations of \cref{alg:scaling_algorithm_nesterov} with $\gamma_k \equiv \frac{e}{L}, \xi_k \equiv \sqrt{\frac{L \Gamma}{\mu e}}$ and $\theta_k \equiv \frac{\sqrt{L\Gamma}}{\sqrt{\mu e} + \sqrt{L\Gamma}}$, we have
\begin{align*}
    \norm{x^{k+1} - x^{*}}_{\hat{D}_{k}}^{2} &\leq \exp\left(-K\sqrt{\frac{\mu e}{4L\Gamma}}\right)\left[\norm{x^{0} - x^{*}}_{\hat{D}_{-1}}^{2} + \frac{2\Gamma}{\mu} (f(x^{0}) - f(x^*))\right].
\end{align*}
\end{theorem}
\end{tcolorbox}
As a result, the next corollary holds.
\begin{tcolorbox}[colback=gray!25, colframe=gray!10, sharp corners]
\begin{corollary}
     Under the conditions of \cref{thm:2}, the required number of iterations $K$ of \cref{alg:scaling_algorithm_nesterov} for reaching $\varepsilon$-accuracy, i.e. $\sqn{x^k - x^*} \leq \varepsilon$, can be upper bounded as
    \begin{align*}
        K = \mathcal{O}\left(\sqrt{\frac{L\Gamma}{\mu e}}\log\left(\frac{\sqn{x_0 - x^*} + \frac{\Gamma}{\mu}(f(x^0) - f(x^*))}{\varepsilon}\right)\right).
    \end{align*}
\end{corollary}
\end{tcolorbox}
\section{Discussion of the results}
Let us highlight the differences between the convergence guarantees of the proposed Algorithms \ref{alg:scaling_algorithm} and \ref{alg:scaling_algorithm_nesterov} compared to their unpreconditioned versions. Two key aspects are worth noting:
\begin{enumerate}
    \item \textbf{The induced norm}. As a convergence criterion in Theorems \ref{thm:1} and \ref{thm:2}, we employ the norm induced by the preconditioning matrix. Essentially, one can transition to the $L_2$-space by applying \cref{prop: scale-bound}, which introduces a multiplicative factor $\Gamma$ related to the initial distance to the optimum. However, this factor is not critical, as it appears inside a \textit{logarithmic} term in the iteration complexity bound on $K$.
    \item \textbf{Exponential term.} In turn, the exponential factor in Theorems \ref{thm:1} and \ref{thm:2} no longer yields a polylogarithmic dependence on $\frac{\Gamma}{e}$ in the final bound on $K$. Specifically, Theorem \ref{thm:1} contributes a term $\frac{\Gamma}{e}$, while Theorem \ref{thm:2} contributes a term $\sqrt{\frac{\Gamma}{e}}$.
\end{enumerate}
As a result, for the heavy-ball method with scaling, no improvement in the dependence on the multiplicative factor is observed -- similar dependence had already been established in prior works \cite{beznosikov2022scaled,sadiev2024stochastic,chezhegov2024local} considering various techniques. This effect can be explained rather straightforwardly: the heavy-ball method does not provide \textit{theoretical} acceleration for the class of functions satisfying Assumptions \ref{asm: strong} and \ref{asm: smooth}. In contrast, the scaled version of Nesterov method demonstrates, to the best of our knowledge, the first improvement of its kind with respect to the factor $\frac{\Gamma}{e}$, as this method enables theoretical acceleration under the given assumptions.

\noindent This phenomenon can be attributed to the fact that the scaling matrix allows one to operate over a modified landscape of the objective function, where, due to \cref{prop: scale-bound}, the effective smoothness and strong convexity constants are altered accordingly; that is, $\mu$ becomes $\frac{\mu}{\Gamma}$, and $L$ becomes $\frac{L}{e}$.
\section{Experiments}
In this section we describe the experimental setups and present the numerical results.
\subsection{Setup}
\textit{Datasets.} We utilize the \texttt{a9a} and \texttt{w8a} LibSVM \cite{chang2011libsvm} datasets in our experiments. These datasets are chosen due to their diverse characteristics and their applicability to classification tasks. We divide the data into training and test parts in a percentage ratio of $80\%$ for training and $20\%$ for testing.

\noindent\textit{Metric. }Since we solve the classification problem, we use standard metrics such as cross-entropy loss. To estimate the rate of convergence we use the square norm of the gradient.

\noindent\textit{Model. }For our experiments, we chose linear model, which in combination with cross-entropy loss handles the binary classification problem effectively.

\noindent\textit{Optimization methods. }For our experiments, we implemented two optimization methods such as Heavy Ball, Nesterov and their scaled versions. Preconditioning matrix is chosen as for the \textsc{Adam} approach. For \cref{alg:scaling_algorithm}, we chose a momentum parameter $\beta_1= 0.9$, for \cref{alg:scaling_algorithm_nesterov} we chose the hyperparameters according to the theory. Learning rate $\gamma$ is chosen as the best option after tuning. 
\subsection{Experiment Results}
\begin{figure}[h]
\label{fig}
    \centering
    \begin{subfigure}{0.49\textwidth}
        \centering
        \includegraphics[width=\textwidth]{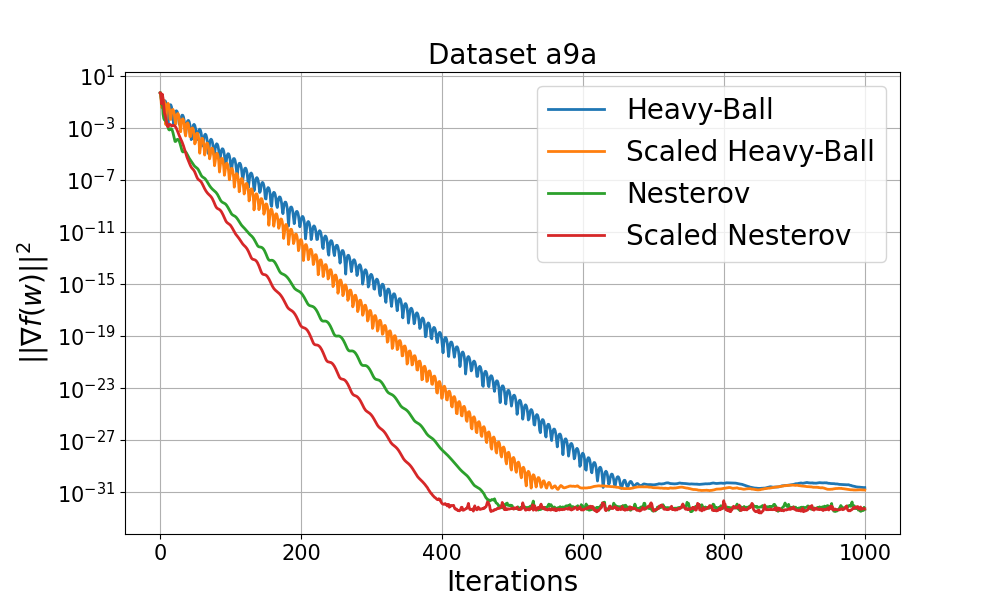}
        \caption{Dataset a9a}
    \end{subfigure}
    \hfill
    \begin{subfigure}{0.49\textwidth}
        \centering
        \includegraphics[width=\textwidth]{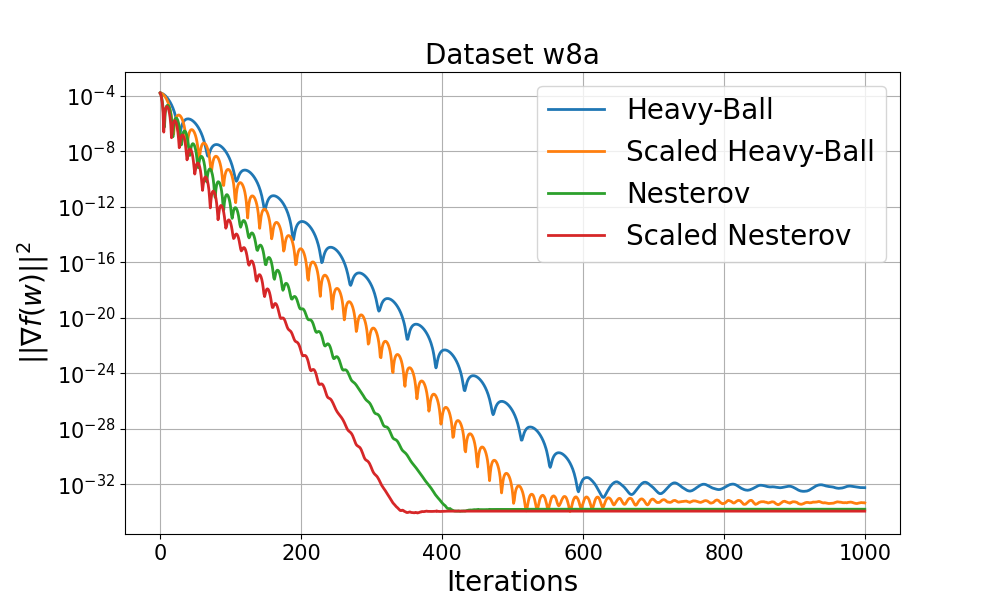}
        \caption{Dataset w8a}
    \end{subfigure}
    \caption{Comparison of Scaled Heavy-Ball and Nesterov methods with unscaled versions.}
\end{figure}
\noindent The results of numerical experiments are presented above. Although the theoretical estimates get worse (due to the $\frac{\Gamma}{e}$ times increase in the condition number of the problem), in practice we see a rather expected effect -- scaled versions of the algorithms allow us to converge faster to the optimum compared to the unpreconditioned techniques. This indicates the applicability of the approach we presented.

\section{Conclusion}
In this work, we proposed the design of two accelerated algorithms incorporating a preconditioning matrix. For the proposed methods \textsc{Preconditioned Heavy Ball} and \textsc{Preconditioned Nesterov}, we also provided theoretical convergence guarantees. These guarantees involve an additional multiplicative factor $\frac{\Gamma}{e}$ that slightly worsens the upper bound on the number of iterations. However, such a factor is standard for methods of this class. Moreover, our empirical results demonstrate that the scaled versions of the algorithms significantly outperform their unscaled counterparts in terms of convergence speed, confirming the practical effectiveness of the proposed approach. The presented study also opens avenues for future research, particularly in the directions of stochastic extensions and analysis under generalized smoothness assumptions.




\bibliographystyle{abbrv}   
\bibliography{ref}

\begin{thebibliography}{10}

\bibitem{beznosikov2022scaled}
A.~Beznosikov, A.~Alanov, D.~Kovalev, M.~Tak{\'a}{\v{c}}, and A.~Gasnikov.
\newblock On scaled methods for saddle point problems.
\newblock {\em arXiv preprint arXiv:2206.08303}, 2022.

\bibitem{cauchy1847methode}
A.~Cauchy et~al.
\newblock M{\'e}thode g{\'e}n{\'e}rale pour la r{\'e}solution des systemes d’{\'e}quations simultan{\'e}es.
\newblock {\em Comp. Rend. Sci. Paris}, 25(1847):536--538, 1847.

\bibitem{chang2011libsvm}
C.-C. Chang and C.-J. Lin.
\newblock Libsvm: a library for support vector machines.
\newblock {\em ACM transactions on intelligent systems and technology (TIST)}, 2(3):1--27, 2011.

\bibitem{chezhegov2024local}
S.~Chezhegov, S.~Skorik, N.~Khachaturov, D.~Shalagin, A.~Avetisyan, M.~Tak{\'a}{\v{c}}, Y.~Kholodov, and A.~Beznosikov.
\newblock Local methods with adaptivity via scaling.
\newblock {\em arXiv preprint arXiv:2406.00846}, 2024.

\bibitem{danilova2022convergence}
M.~Danilova.
\newblock On the convergence analysis of aggregated heavy-ball method.
\newblock In {\em International Conference on Mathematical Optimization Theory and Operations Research}, pages 3--17. Springer, 2022.

\bibitem{dozat2016incorporating}
T.~Dozat.
\newblock Incorporating nesterov momentum into adam.
\newblock 2016.

\bibitem{duchi2011}
J.~Duchi, E.~Hazan, and Y.~Singer.
\newblock Adaptive subgradient methods for online learning and stochastic optimization.
\newblock {\em Journal of Machine Learning Research}, 12:2121--2159, 2011.

\bibitem{flammarion2015averaging}
N.~Flammarion and F.~Bach.
\newblock From averaging to acceleration, there is only a step-size.
\newblock In {\em Conference on learning theory}, pages 658--695. PMLR, 2015.

\bibitem{fletcher1963rapidly}
R.~Fletcher and M.~J. Powell.
\newblock A rapidly convergent descent method for minimization.
\newblock {\em The computer journal}, 6(2):163--168, 1963.

\bibitem{ghadimi2012optimal}
S.~Ghadimi and G.~Lan.
\newblock Optimal stochastic approximation algorithms for strongly convex stochastic composite optimization i: A generic algorithmic framework.
\newblock {\em SIAM Journal on Optimization}, 22(4):1469--1492, 2012.

\bibitem{goldfarb1970family}
D.~Goldfarb.
\newblock A family of variable-metric methods derived by variational means.
\newblock {\em Mathematics of computation}, 24(109):23--26, 1970.

\bibitem{goodfellow2016deep}
I.~Goodfellow, Y.~Bengio, and A.~Courville.
\newblock Deep feedforward networks.
\newblock {\em Deep learning}, 1:161--217, 2016.

\bibitem{jahani2021doubly}
M.~Jahani, S.~Rusakov, Z.~Shi, P.~Richt{\'a}rik, M.~W. Mahoney, and M.~Tak{\'a}{\v{c}}.
\newblock Doubly adaptive scaled algorithm for machine learning using second-order information.
\newblock {\em arXiv preprint arXiv:2109.05198}, 2021.

\bibitem{kingma2014}
D.~P. Kingma and J.~Ba.
\newblock Adam: A method for stochastic optimization.
\newblock {\em arXiv preprint arXiv:1412.6980}, 2014.

\bibitem{liu1989limited}
D.~C. Liu and J.~Nocedal.
\newblock On the limited memory bfgs method for large scale optimization.
\newblock {\em Mathematical programming}, 45(1):503--528, 1989.

\bibitem{liu2023sophia}
H.~Liu, Z.~Li, D.~Hall, P.~Liang, and T.~Ma.
\newblock Sophia: A scalable stochastic second-order optimizer for language model pre-training.
\newblock {\em arXiv preprint arXiv:2305.14342}, 2023.

\bibitem{loshchilov2017decoupled}
I.~Loshchilov.
\newblock Decoupled weight decay regularization.
\newblock {\em arXiv preprint arXiv:1711.05101}, 2017.

\bibitem{nesterov2013introductory}
Y.~Nesterov.
\newblock {\em Introductory lectures on convex optimization: A basic course}, volume~87.
\newblock Springer Science \& Business Media, 2013.

\bibitem{nesterov1983}
Y.~E. Nesterov.
\newblock A method for unconstrained convex minimization problem with the rate of convergence \(o(1/k^2)\).
\newblock {\em Doklady AN SSSR}, 269:543--547, 1983.

\bibitem{polyak1964}
B.~T. Polyak.
\newblock Some methods of speeding up the convergence of iteration methods.
\newblock {\em USSR Computational Mathematics and Mathematical Physics}, 4(5):1--17, 1964.

\bibitem{reddi2020adaptive}
S.~Reddi, Z.~Charles, M.~Zaheer, Z.~Garrett, K.~Rush, J.~Kone{\v{c}}n{\`y}, S.~Kumar, and H.~B. McMahan.
\newblock Adaptive federated optimization.
\newblock {\em arXiv preprint arXiv:2003.00295}, 2020.

\bibitem{reddi2019convergence}
S.~J. Reddi, S.~Kale, and S.~Kumar.
\newblock On the convergence of adam and beyond.
\newblock {\em arXiv preprint arXiv:1904.09237}, 2019.

\bibitem{sadiev2024stochastic}
A.~Sadiev, A.~Beznosikov, A.~J. Almansoori, D.~Kamzolov, R.~Tappenden, and M.~Tak{\'a}{\v{c}}.
\newblock Stochastic gradient methods with preconditioned updates.
\newblock {\em Journal of Optimization Theory and Applications}, 201(2):471--489, 2024.

\bibitem{shalev2014understanding}
S.~Shalev-Shwartz and S.~Ben-David.
\newblock {\em Understanding machine learning: From theory to algorithms}.
\newblock Cambridge university press, 2014.

\bibitem{shanno1970conditioning}
D.~F. Shanno.
\newblock Conditioning of quasi-newton methods for function minimization.
\newblock {\em Mathematics of computation}, 24(111):647--656, 1970.

\bibitem{tieleman2012lecture}
T.~Tieleman.
\newblock Lecture 6.5-rmsprop: Divide the gradient by a running average of its recent magnitude.
\newblock {\em COURSERA: Neural networks for machine learning}, 4(2):26, 2012.

\bibitem{yao2021adahessian}
Z.~Yao, A.~Gholami, S.~Shen, M.~Mustafa, K.~Keutzer, and M.~Mahoney.
\newblock Adahessian: An adaptive second order optimizer for machine learning.
\newblock In {\em proceedings of the AAAI conference on artificial intelligence}, volume~35, pages 10665--10673, 2021.

\end{thebibliography}
\newpage
\appendix
\section{Supplementary materials}
Before we start, let us formulate auxiliary propositions.
\begin{proposition}[Young's inequality]
\label{prop: young}
For all $x, y \in \mathbb{R}^d$ and for any $\lambda > 0$ the next inequality holds:
\begin{align*}
    \ev{x, y} \leq \frac{\norm{x}^2}{2\lambda} + \frac{\lambda\norm{y}^2}{2}.
\end{align*}
\end{proposition}
\begin{proposition}[Jensen's inequality]
\label{prop: jensen}
For any convex $f: \mathbb{R}^d \rightarrow \mathbb{R}$ and $\{a_i\}_{i=1}^n: a_i \geq 0$ and $\sum\limits_{i=1}^n a_i = 1$, the next inequality holds:
\begin{align*}
    f\left(\sum\limits_{i=1}^na_i x_i \right) \leq \sum\limits_{i=1}^na_if\left(x_i \right).
\end{align*}    
\end{proposition}
\subsection{Proof of Preconditioned Heavy-Ball method}
\label{hb}
\begin{tcolorbox}[colback=gray!25, colframe=gray!10, sharp corners]
\begin{lemma}[\textbf{Descent lemma}] 
\label{lem: descent-hb}
Suppose that Assumptions \ref{asm: strong}, \ref{asm: smooth} and \ref{asm: scaling-matr} hold. Then, if the sequence $\{x_k\}_{k=0}$ is produced by \cref{alg:scaling_algorithm} and $\gamma$ satisfies
\begin{align*}
F := \frac{\gamma}{1 - \beta_1} \leq \frac{e}{4L},
\end{align*} 
the next inequality holds:
\begin{align*}
\frac{F}{2} \left( f(x_k) - f(x^*) \right) &\leq \left( 1 - \frac{\mu F}{4\Gamma} \right) \norm{\tilde{x}_k - x^*}_{\hat{D}_{k-1}}^2 - \norm{\tilde{x}_{k+1} - x^*}_{\hat{D}_{k}}^2 \\&+ \frac{3LF}{e} \norm{x_k - \tilde{x}_k}_{\hat{D}_{k}}^2
\end{align*}
for all \( k \geq 0 \), where $\tilde{x}_k = x_k - \frac{\beta_1\gamma}{1 - \beta_1} V_{k-1}$.
\end{lemma}
\end{tcolorbox}
\begin{proof}
Starting with the notation of virtual sequence $\{\tilde{x}_k\}$ and the update rule:
\begin{align}
\label{eq: virtual-hb}
     \tilde{x}_{k+1} &= x_{k+1} - \frac{\beta_1\gamma}{1 - \beta_1}\cdot V_{k} = x_{k} - \gamma V_{k} - \frac{\beta_1\gamma}{1-\beta_1}V_{k} = x_{k} - \frac{(1 - \beta_1)\gamma + \beta_1\gamma}{1 - \beta_1}V_{k} \nonumber\\&= x_{k} - \frac{\gamma}{1 - \beta_1}V_{k} = x_{k} - \frac{\gamma}{1 - \beta_1}\Big(\beta_1 V_{k-1} + \hat{D}_{k}^{-1}\nabla f(x_{k})\Big) \nonumber \\
     &= x_{k} - \frac{\gamma\beta_1}{1 - \beta_1}V_{k-1} - \frac{\gamma}{1 - \beta_1}\hat{D}_{k}^{-1}\nabla f(x_{k}) = \tilde{x}_{k} - \frac{\gamma}{1 - \beta_1}\hat{D}_{k}^{-1}\nabla f(x_{k}).
\end{align}
Therefore, using \eqref{eq: virtual-hb}, we get
\begin{align}
\label{eq: hb-descent}
    \norm{\tilde{x}_{k+1}  - x^{*}}_{\hat{D}_{k}}^{2} &= \norm{\tilde{x}_{k} - x^{*}}_{\hat{D}_{k}}^{2} + F^{2}\norm{\hat{D}_{k}^{-1}\nabla f(x_{k})}_{\hat{D}_{k}}^{2} \nonumber\\&- 2F\ev{\tilde{x}_{k} - x^{*}, \hat{D}_{k}\hat{D}_{k}^{-1}\nabla f(x_{k})} \nonumber\\&= \norm{\tilde{x}_{k} - x^{*}}_{\hat{D}_{k}}^{2} + F^{2}\norm{\nabla f(x_{k})}_{\hat{D}_{k}^{-1}}^{2} - 2F\ev{{x}_{k} - x^{*}, \nabla f(x_{k})} \nonumber\\&
    - 2F\ev{\tilde{x}_{k} - x_{k}, \nabla f(x_{k})}.
\end{align}
To upper bound the second term in \eqref{eq: hb-descent}, we apply \cref{prop: equiv-L} and \cref{prop: scale-bound}:
\begin{align}
\label{eq: bound-norm-grad}
    \norm{\nabla f(x_{k})}_{\hat{D}_{k}^{-1}}^{2}
    \leq \frac{1}{e}\sqn{\nabla f(x_{k})}_2 \leq \frac{2L}{e}\left(f(x_{k}) - f(x^{*})\right).
\end{align}
For the third term in \eqref{eq: hb-descent}, according to \cref{asm: strong}, one can obtain
\begin{align}
\label{eq: hb-first-ev}
    \ev{x_{k} - x^{*}, \nabla f(x_{k})} &\geq f(x_{k}) - f(x^{*}) + \frac{\mu}{2}\norm{x_{k} - x^{*}}^{2}_2 \nonumber\\&\geq f(x_{k}) - f(x^{*}) + \frac{\mu}{2\Gamma}\norm{x_{k} - x^{*}}_{\hat{D}_k}^{2}.
\end{align}
For the last term in \eqref{eq: hb-descent}, we use the Young's inequality (\cref{prop: young}) with $\lambda = \frac{1}{2L}$ and \cref{prop: equiv-L}:
\begin{align}
\label{eq: hb-second-ev}
    -2F\ev{\tilde{x}_{k} - x_{k}, \nabla f(x_{k})} &\leq 2LF\norm{\tilde{x}_{k} - x_{k}}_{2}^{2} + \frac{F}{2L}\norm{\nabla f(x_{k})}_{2}^{2}
    \nonumber\\&\leq \frac{2LF}{e}\norm{\tilde{x}_{k} - x_{k}}_{\hat{D}_k}^{2} + F\left(f(x_{k}) - f(x^{*})\right).
\end{align}
Substituting \eqref{eq: bound-norm-grad}, \eqref{eq: hb-first-ev} and \eqref{eq: hb-second-ev} into \eqref{eq: hb-descent}, we obtain 
\begin{align*}
    \norm{\tilde{x}_{k+1}  - x^{*}}_{\hat{D}_{k}}^{2} &= \norm{\tilde{x}_{k} - x^{*}}_{\hat{D}_{k}}^{2} + F^{2}\norm{\nabla f(x_{k})}_{\hat{D}_{k}^{-1}}^{2} - 2F\ev{\tilde{x}_{k} - x^{*}, \nabla f(x_{k})}
    \\&- 2F\ev{\tilde{x}_{k} - x_{k}, \nabla f(x_{k})} \\
    &\leq \norm{\tilde{x}_{k} - x^{*}}_{\hat{D}_{k}}^{2} + \frac{2LF^2}{e}\left(f(x_{k}) - f(x_{*})\right) \\&- 2F\left(f(x_{k}) - f(x^{*}) + \frac{\mu}{2\Gamma}\norm{x_{k} - x^{*}}_{\hat{D}_{k}}^{2}\right) \\
    &+ \frac{2LF}{e}\norm{\tilde{x}_{k} - x_{k}}_{\hat{D}_k}^{2} + F\left(f(x_{k}) - f(x^{*})\right) \\&= \norm{\tilde{x}_{k} - x^{*}}_{\hat{D}_{k}}^{2} + \left(\frac{2LF^{2}}{e} - F\right)\left(f(x_{k}) - f(x^{*})\right)\\&+ \frac{2LF}{e}\norm{\tilde{x}_{k} - x_{k}}_{\hat{D}_k}^{2} - \frac{\mu F}{\Gamma}\norm{x_k - x^*}^2_{\hat{D}_k}. \\
\end{align*}
To bound the last term, we apply one-dimensional Young's inequality (\cref{prop: young}) with $\lambda = 1$ to obtain
\begin{align*}
    \frac{\mu F}{2 \Gamma}\norm{\tilde{x}_k - x^*}^2_{\hat{D}_k} &= \frac{\mu F}{2 \Gamma}\sum\limits_{i=1}^d [\hat{D}_k]_{ii}[\tilde{x}_k - x^*]_i^2 \\&\leq \frac{\mu F}{2 \Gamma}\sum\limits_{i=1}^d [\hat{D}_k]_{ii}\left(2[x_k - x^*]_i^2 + 2[\tilde{x}_k - x_k]_i^2\right) \\&= \frac{\mu F}{\Gamma}\norm{\tilde{x}_k - x_k}_{\hat{D}_k}^2 + \frac{\mu F}{\Gamma}\norm{x_k - x^*}_{\hat{D}_k}^2. 
\end{align*}
Therefore, one can get
\begin{align*}
    \norm{\tilde{x}_{k+1}  - x^{*}}_{\hat{D}_{k}}^{2} &\leq \left(1 - \frac{\mu F}{2\Gamma}\right)\norm{\tilde{x}_{k} - x^{*}}_{\hat{D}_{k}}^{2} + \left(\frac{2LF^{2}}{e} - F\right)\left(f(x_{k}) - f(x^{*})\right) \\&+ \left(\frac{2LF}{e} + \frac{\mu F}{\Gamma}\right)\norm{\tilde{x}_{k} - x_{k}}_{\hat{D}_k}^{2} \\&\leq \left(1 - \frac{\mu F}{2\Gamma}\right)\norm{\tilde{x}_{k} - x^{*}}_{\hat{D}_{k}}^{2} + \left(\frac{2LF^{2}}{e} - F\right)\left(f(x_{k}) - f(x^{*})\right) \\&+ \frac{3LF}{e}\norm{\tilde{x}_{k} - x_{k}}_{\hat{D}_k}^{2},
\end{align*}
where in the last inequality we applied $\frac{\mu}{\Gamma} \leq \frac{L}{e}$. Due to the bound on $F$ we have $\frac{2LF^2}{e} \leq \frac{F}{2}$. What is more, applying \cref{prop: ultra} with  $\beta_2 \geq 1 - \frac{\mu F}{4\Gamma C}$ allows to get
\begin{align*}
    \norm{\tilde{x}_{k} - x^{*}}_{\hat{D}_{k}}^{2} \leq \left(1 + \frac{\mu F}{4\Gamma}\right)\norm{\tilde{x}_{k} - x^{*}}_{\hat{D}_{k-1}}^{2}. 
\end{align*}
Combining this with $(1 - x)\left(1 + \frac{x}{2}\right) \leq 1 - \frac{x}{2}$, one can obtain
\begin{align*}
     \norm{\tilde{x}_{k+1}  - x^{*}}_{\hat{D}_{k}}^{2} &\leq \left(1 - \frac{\mu F}{4\Gamma}\right)\norm{\tilde{x}_{k} - x^{*}}_{\hat{D}_{k-1}}^{2} - \frac{F}{2}\left(f(x_{k}) - f(x^{*})\right) \\&+ \frac{3LF}{e}\norm{\tilde{x}_{k} - x_{k}}_{\hat{D}_k}^{2},
\end{align*}
what concludes the proof.  
\end{proof}
\begin{tcolorbox}[colback=gray!25, colframe=gray!10, sharp corners]
\begin{lemma}[Auxiliary lemma] 
\label{lem: aux-hb}
Suppose that Assumptions \ref{asm: smooth} and \ref{asm: scaling-matr} hold. Thus, if the sequence $\{x_k\}_{k=0}$ is produced by \cref{alg:scaling_algorithm}, the following inequality is satisfied for all $k$:
\begin{align*}
    \norm{\tilde{x}_{k} - x_{k}}_{\hat{D}_k}^{2} \leq \frac{2L\beta_1^2\gamma^2}{(1 - \beta_1)^3e} \sum\limits_{t=0}^{k-1} \beta_1^{k - 1 - t} \left(1 + \frac{\mu F}{4\Gamma}\right)^{k - t} (f(x_t) - f(x^*)). 
\end{align*}
\end{lemma}
\end{tcolorbox}
\begin{proof}
Let us substitute the analytical form of $\tilde{x_{k}}$:
\begin{align}
\label{eq: aux-start}
    \norm{\tilde{x_k} - x_k}^2_{\hat{D}_{k}} &= \norm{\frac{\beta_1\gamma}{1 - \beta_1}V_{k-1}}^2_{\hat{D}_{k}} = \norm{\frac{\beta_1\gamma}{1 - \beta_1}\sum\limits_{t=0}^{k-1}\beta_1^{k - 1 - t}\hat{D}_t^{-1}\nabla f(x_t)}^2_{\hat{D}_{k}} \nonumber\\&= \frac{\beta_1^2\gamma^2}{(1 - \beta_1)^2}\left(\sum\limits_{j=0}^{k-1}\beta_1^{k - 1 - j}\right)^2\norm{\sum\limits_{t=0}^{k-1}\frac{\beta_1^{k - 1 - t}}{\sum\limits_{j=0}^{k-1}\beta_1^{k - 1 - j}}\hat{D}_t^{-1}\nabla f(x_t)}^2_{\hat{D}_{k}} \nonumber\\&\leq \frac{\beta_1^2\gamma^2}{(1 - \beta_1)^2}\left(\sum\limits_{j=0}^{k-1}\beta_1^{k - 1 - j}\right)\sum\limits_{t=0}^{k-1}\beta_1^{k-1-t}\norm{\hat{D}_t^{-1}\nabla f(x_t)}^2_{\hat{D}_{k}}, 
\end{align}
where we apply the Jensen's inequality (\cref{prop: jensen}) for the convex function $\|\cdot\|^2$. Moreover, due to \cref{prop: ultra} with $\beta_2 \geq 1 - \frac{\mu F}{4\Gamma C}$, we have
\begin{align}
\label{eq: bound t-k}
    \norm{\hat{D}_t^{-1}\nabla f(x_t)}^2_{\hat{D}_{k}} &\leq \left(1 + \frac{\mu F}{4\Gamma}\right) \norm{\hat{D}_t^{-1}\nabla f(x_t)}^2_{\hat{D}_{k-1}} \leq \ldots \nonumber\\&\leq \left(1 + \frac{\mu F}{4\Gamma}\right)^{k - t}\norm{\hat{D}_t^{-1}\nabla f(x_t)}^2_{\hat{D}_{t}} \nonumber\\&= \left(1 + \frac{\mu F}{4\Gamma}\right)^{k - t}\norm{\nabla f(x_t)}^2_{\hat{D}_t^{-1}} \nonumber\\&\leq \left(1 + \frac{\mu F}{4\Gamma}\right)^{k - t}\frac{1}{e}\norm{\nabla f(x_t)}^2_{2} \nonumber\\&\leq \left(1 + \frac{\mu F}{4\Gamma}\right)^{k - t}\frac{2L}{e}\left(f(x_t) - f(x^*)\right),  
\end{align}
where in the last inequality we applied \cref{prop: equiv-L} and \cref{prop: scale-bound}. Substituting \eqref{eq: bound t-k} into \eqref{eq: aux-start}, we finish the proof.  
\end{proof}
\noindent Now we are ready to formulate the theorem of convergence of Preconditioned Heavy-Ball.
\begin{tcolorbox}[colback=gray!25, colframe=gray!10, sharp corners]
\begin{theorem}
Suppose that Assumptions \ref{asm: strong}, \ref{asm: smooth} and \ref{asm: scaling-matr}  hold. Then, after $K$ iterations of \cref{alg:scaling_algorithm} with $\gamma = \frac{(1 - \beta_1)^2e}{12 L}$, we have
\begin{align*}
     f\left(\frac{1}{W_{K-1}}\sum\limits_{k=0}^{K-1} w_kx_k\right) - f(x^*) &\leq 4\exp\left(-\frac{(1 - \beta_1)\mu e K}{48 L \Gamma}\right) L\norm{{x}_{0} - x^{*}}_{\hat{D}_{-1}}^{2},
\end{align*}
where $W_{K-1} = \sum_{k=0}^{K-1} w_k$, $w_k = \left(1 - \frac{\mu F}{4 \Gamma}\right)^{-(k+1)}$.
\end{theorem}
\end{tcolorbox}
\begin{proof}
Let us denote $w_k = \left(1 - \frac{\mu F}{4 \Gamma}\right)^{-(k+1)}$. Therefore, summing the results from \cref{lem: descent-hb} with weights $w_k$, we get
\begin{align*}
    \sum\limits_{k=0}^{K-1} &\frac{w_k F}{2} (f(x_k) - f(x^*)) \\&\leq \sum\limits_{k=0}^{K-1} w_{k}\Bigg[\left(\left(1 - \frac{\mu F}{4\Gamma}\right)\norm{\tilde{x}_{k} - x^{*}}_{\hat{D}_{k-1}}^{2} - \norm{\tilde{x}_{k+1}  - x^{*}}_{\hat{D}_{k}}^{2}\right) \\&+ \frac{3LF}{e}\norm{\tilde{x}_{k} - x_{k}}_{\hat{D}_k}^{2} \Bigg]\\&=
    \sum\limits_{k=0}^{K-1} \Bigg[w_{k-1}\norm{\tilde{x}_{k} - x^{*}}_{\hat{D}_{k-1}}^{2} - w_{k}\norm{\tilde{x}_{k+1}  - x^{*}}_{\hat{D}_{k}}^{2} + \frac{3LFw_k}{e}\norm{\tilde{x}_{k} - x_{k}}_{\hat{D}_k}^{2}\Bigg].
\end{align*}
The main question arises for the last term. Applying \cref{lem: aux-hb}, one can obtain
\begin{align*}
    \sum\limits_{k=0}^{K-1} &\frac{3LFw_k}{e}\norm{\tilde{x}_{k} - x_{k}}_{\hat{D}_k}^{2} \\&\leq \frac{6L^2F \beta_1^2\gamma^2}{(1 - \beta_1)^3e^2}\sum\limits_{k=0}^{K-1} \sum\limits_{t=0}^{k-1} w_k \beta_1^{k - 1 - t} \left(1 + \frac{\mu F}{4\Gamma}\right)^{k - t} (f(x_t) - f(x^*)).  
\end{align*}
Now we decompose $w_k$ as
\begin{align*}
    w_k = w_t\left(1 - \frac{\mu F}{4\Gamma}\right)^{- {(k - t)}} \leq w_t\left(1 + \frac{\mu F}{2\Gamma}\right)^{{k - t}} \leq w_t\left(1 + \frac{1 - \beta_1}{2}\right)^{k-t}, 
\end{align*}
where the last inequality holds due to the choice of $\gamma$. Consequently, using that $\beta_1 = 1 - (1 - \beta_1)$ and $(1 - x)\left(1 + \frac{x}{2}\right) \leq \left(1  - \frac{x}{2}\right)$, we have
\begin{align*}
    \sum\limits_{k=0}^{K-1} &\frac{3LFw_k}{e}\norm{\tilde{x}_{k} - x_{k}}_{\hat{D}_k}^{2} \\&\leq \frac{6L^2F \beta_1\gamma^2}{(1 - \beta_1)^3e^2}\sum\limits_{k=0}^{K-1} \sum\limits_{t=0}^{k-1} w_t \left(1 - \frac{1 -\beta_1}{2}\right)^{k - t} \left(1 + \frac{\mu F}{4\Gamma}\right)^{k - t} (f(x_t) - f(x^*)).
\end{align*}
Moreover, since $1 - \beta_1 \geq \frac{\mu F}{\Gamma}$, we get
\begin{align*}
    \sum\limits_{k=0}^{K-1} &\frac{3LFw_k}{e}\norm{\tilde{x}_{k} - x_{k}}_{\hat{D}_k}^{2} \\&\leq \frac{6L^2F \beta_1\gamma^2}{(1 - \beta_1)^3e^2}\sum\limits_{k=0}^{K-1} \sum\limits_{t=0}^{k-1} w_t \left(1 - \frac{1 -\beta_1}{4}\right)^{k - t} (f(x_t) - f(x^*)).
\end{align*}
Therefore, we obtain
\begin{align*}
    &\sum\limits_{k=0}^{K-1} \frac{w_k F}{2}(f(x_k) - f(x^*)) \\&\leq
    w_{-1}\norm{\tilde{x}_{0} - x^{*}}_{\hat{D}_{-1}}^{2} + \frac{6L^2F \beta_1\gamma^2}{(1 - \beta_1)^3e^2}\sum\limits_{k=0}^{K-1} \sum\limits_{t=0}^{k-1} w_t \left(1 - \frac{1 -\beta_1}{4}\right)^{k - t} (f(x_t) - f(x^*)).
\end{align*}
It can be easily shown that the coefficient related to $(f(x_r) - f(x^*))$ in the right-hand side can be upper bounded as
\begin{align*}
    \frac{6L^2 F\beta_1 \gamma^2}{(1 - \beta_1)^3 e^2}\sum\limits_{k=0}^{\infty} \left(1 - \frac{1 -\beta_1}{4}\right)^{k} \leq \frac{24L^2 F\beta_1 \gamma^2}{(1 - \beta_1)^4 e}w_r = \frac{24L^2 F^3\beta_1}{(1 - \beta_1)^2 e^2}w_r.
\end{align*}
Applying $F \leq \frac{(1 - \beta_1)e}{12L}$, we obtain
\begin{align*}
     \frac{24L^2 F^3\beta_1}{(1 - \beta_1)^2 e^2} \leq \frac{F}{4}.
\end{align*}
As a result, one can get
\begin{align*}
    \sum\limits_{k=0}^{K-1} \frac{w_k F}{4} (f(x_k) - f(x^*)) &\leq
    w_{-1}\norm{\tilde{x}_{0} - x^{*}}_{\hat{D}_{-1}}^{2}.
\end{align*}
Dividing both sides by $W_{K-1} = \sum_{k=0}^{K-1} w_k$ and using that $\sum_{k=0}^{K-1} w_k \geq w_{K-1} \geq \exp\left(\frac{\mu F K}{4 \Gamma}\right)$, we have.
\begin{align*}
    \frac{1}{W_{K-1}}\sum\limits_{k=0}^{K-1} w_k (f(x_k) - f(x^*)) &\leq 4\exp\left(-\frac{\mu F K}{4 \Gamma}\right) \norm{\tilde{x}_{0} - x^{*}}_{\hat{D}_{-1}}^{2}.
\end{align*}
Applying the Jensen's inequality (\cref{prop: jensen}) to the left-hand side and substituting the choice of $\gamma$, we conclude the proof.  
\end{proof}

\subsection{Proof of Preconditioned Nesterov Method}
\label{nesterov}
\begin{tcolorbox}[colback=gray!25, colframe=gray!10, sharp corners]
\begin{lemma}[Auxiliary lemma]
\label{lem:lemma1}
Suppose that Assumptions \ref{asm: strong}, \ref{asm: smooth} and \ref{asm: scaling-matr}  hold. Then, if the sequence $\{x_k\}_{k=0}$ is produced by \cref{alg:scaling_algorithm_nesterov}, for all $k$ and any $u \in \mathbb{R}^d$ the following inequality is satisfied:
\begin{align*}
f(x_f^{k + 1}) &\leq f(u) - \ev{ \nabla f(x_{g}^{k}), u - x_{g}^{k} } - \frac{\mu}{2} \norm{u - x_{g}^{k}}_{2}^{2} \\&+ \left(\frac{L\gamma_k}{2e} - 1\right)\gamma_k\norm{\nabla f(x_g^k)}^2_{\hat{D}_{k}^{-1}}.
\end{align*}
\end{lemma}
\end{tcolorbox}
\begin{proof}
Start with the $L$-smoothness of the function $f$:
\begin{align*}
f(x_f^{k + 1}) \leq f(x_{g}^{k})+\langle \nabla f(x_{g}^{k}), x_{f}^{k + 1} - x_{g}^{k}\rangle + \frac{L}{2}\norm{x_{f}^{k + 1} - x_{g}^{k}}_{2}^{2}. 
\end{align*}
Applying \cref{prop: scale-bound} to the last term, we obtain
\begin{align}
\label{eq: start-aux-1}
    f(x_f^{k + 1}) \leq f(x_{g}^{k}) + \ev{\nabla f(x_{g}^{k}), x_{f}^{k + 1} -x_{g}^{k}} + \frac{L}{2e}\norm{x_{f}^{k + 1} -x_{g}^{k}}_{\hat{D}_{k}}^{2}.
\end{align}
The last term can be decomposed due to the update rule of the algorithm:
\begin{align}
\label{eq: aux-1}
    \norm{x_{f}^{k + 1} -x_{g}^{k}}_{\hat{D}_{k}}^{2} = \gamma_k^2\norm{\nabla f(x_g^k)}^2_{\hat{D}_{k}^{-1}}.
\end{align}
Moreover, the same update rule can be applied to the second term:
\begin{align}
\label{eq: aux-2}
    \ev{\nabla f(x_{g}^{k}), x_{f}^{k + 1} -x_{g}^{k}} = -\gamma_k\norm{\nabla f(x_g^k)}^2_{\hat{D}_{k}^{-1}}.
\end{align}
The first term can be upper bounded with the \cref{asm: strong}:
\begin{align}
\label{eq: aux-3}
    f(x_{g}^{k}) \leq f(u) - \ev{ \nabla f(x_{g}^{k}), u - x_{g}^{k} } - \frac{\mu}{2} \norm{u - x_{g}^{k}}_{2}^{2}.
\end{align}
Substituting \eqref{eq: aux-1}, \eqref{eq: aux-2} and \eqref{eq: aux-3} into \eqref{eq: aux-start} gives the final result. 
\end{proof}

\begin{tcolorbox}[colback=gray!25, colframe=gray!10, sharp corners]
\begin{lemma}[\textbf{Descent lemma}]
\label{lem: descent-nesterov}
Suppose that Assumptions \ref{asm: strong}, \ref{asm: smooth} and \ref{asm: scaling-matr}  hold. Then, if the sequence $\{x_k\}_{k=0}$ is produced by \cref{alg:scaling_algorithm_nesterov}, for all $k$ with $\xi_k \geq 1, \frac{\xi_k^2\gamma_k\mu}{\Gamma} \geq 1$, $\gamma_k \leq \frac{e}{L}$ and $\theta_k = \frac{\xi_k}{1 + \xi_k}$, the next inequality is satisfied:
\begin{align*}
    \norm{x^{k+1} - x^{*}}_{\hat{D}_{k}}^{2} &+ 2\gamma_k\xi_k^2 (f(x_f^{k+1}) - f(x^*)) \\&\leq \left(1 - \frac{1}{\xi_k}\right)\left[\norm{x^{k} - x^{*}}_{\hat{D}_{k}}^{2} + 2\gamma_k\xi_k^2 (f(x_f^{k}) - f(x^*))\right].
\end{align*}
\end{lemma}
\end{tcolorbox}
\begin{proof} 
Starting with the update rule:
\begin{align}
\label{eq: descent-nest-start}
\norm{x^{k+1} - x^{*}}_{\hat{D}_{k}}^{2} &= \norm{\xi_{k}x_{f}^{k+1} + (1 - \xi_{k})x_{f}^{k} - x^{*}}_{\hat{D}_{k}}^{2} \nonumber\\&= \norm{\xi_{k}\left(x_{g}^{k} - \gamma_{k}\hat{D}_{k}^{-1}\nabla f(x_{g}^{k})\right) + (1 - \xi_{k})x_{f}^{k} - x^{*}}_{\hat{D}_{k}}^{2} \nonumber\\
&= \norm{\xi_{k}x_{g}^{k} + (1 - \xi_{k})x_{f}^{k} - x^{*}}_{\hat{D}_{k}}^{2} + \gamma_{k}^{2}\xi_{k}^{2}\norm{\hat{D}_{k}^{-1}\nabla f(x_{g}^{k})}_{\hat{D}_{k}}^{2}\nonumber\\
&-2\gamma_{k}\xi_{k}\ev{\hat{D}_{k}^{-1}\nabla f(x_{g}^{k}), \xi_{k}x_{g}^{k} + (1 - \xi_{k})x_{f}^{k} - x^{*}}_{\hat{D}_{k}} \nonumber\\&= \norm{\xi_{k}x_{g}^{k} + (1 - \xi_{k})x_{f}^{k} - x^{*}}_{\hat{D}_{k}}^{2} + \gamma_{k}^{2}\xi_{k}^{2}\norm{\nabla f(x_{g}^{k})}_{\hat{D}_{k}^{-1}}^{2}\nonumber\\
&-2\gamma_{k}\xi_{k}\ev{\nabla f(x_{g}^{k}), \xi_{k}x_{g}^{k} + (1 - \xi_{k})x_{f}^{k} - x^{*}}.
\end{align}
Let us decompose the first term. According to the update rule, we get
\begin{align*}
    \xi_{k}x_{g}^{k} + (1 - \xi_{k})x_{f}^{k} - x^{*} &= \xi_{k}x_{g}^{k} + \frac{(1 - \xi_{k})}{\theta_k}\theta_k x_{f}^{k} - x^{*} \\&= \xi_{k}x_{g}^{k} + \frac{(1 - \xi_{k})}{\theta_k}(x_g^k - (1 - \theta_k)x^k) - x^{*}.
\end{align*}
Choosing $\theta_k$ as $\frac{\xi_k}{1 + \xi_k}$, after simple estimations, one can obtain
\begin{align*}
    \xi_{k}x_{g}^{k} + \frac{(1 - \xi_{k})}{\theta_k}(x_g^k - (1 - \theta_k)x^k) - x^{*} = x^k - \frac{1}{\xi_k}(x^k - x_g^k).
\end{align*}
Hence, the first term in \eqref{eq: descent-nest-start} can be decomposed as
\begin{align}
\label{eq: nest-1}
    \norm{\xi_{k}x_{g}^{k} + (1 - \xi_{k})x_{f}^{k} - x^{*}}_{\hat{D}_{k}}^{2} &= \norm{x^k - \frac{1}{\xi_k}(x^k - x_g^k) - x^{*}}_{\hat{D}_{k}}^{2} \nonumber\\&= \norm{x^k - x^{*}}^2_{\hat{D}_k} - \frac{2}{\xi_k}\ev{x^k - x^*, x^k - x_g^k}_{\hat{D}_k} \nonumber\\&+ \frac{1}{\xi_k^2}\norm{x^k - x_g^k}^2_{\hat{D}_k}. 
\end{align}
Let us note that
\begin{align*}
    -2\ev{x^k - x^*, x_g^k- x^k}_{\hat{D}_k} = \norm{x_g^k - x^k}^2_{\hat{D}_k} + \norm{x^k - x^*}^2_{\hat{D}_k} - \norm{x_g^k - x^*}^2_{\hat{D}_k}.
\end{align*}
Thus, continue with \eqref{eq: nest-1}:
\begin{align}
\label{eq: n-first}
    \norm{\xi_{k}x_{g}^{k} + (1 - \xi_{k})x_{f}^{k} - x^{*}}_{\hat{D}_{k}}^{2} &= \left(1 - \frac{1}{\xi_k}\right)\norm{x^k - x^*}^2_{\hat{D}_k} \nonumber\\&+ \left(\frac{1}{\xi^2} -\frac{1}{\xi}\right)\norm{x^k - x_g^k}^2_{\hat{D}_k} + \frac{1}{\xi_k}\norm{x_g^k - x^*}^2_{\hat{D}_k}.
\end{align}
From \cref{lem:lemma1} with $u = x_f^k$ and $u = x^*$, under \cref{prop: scale-bound} we get
\begin{align}
\label{eq: n-sec}
    f(x_f^{k + 1}) &\leq f(x_f^k) - \ev{ \nabla f(x_{g}^{k}), x_f^k - x_{g}^{k} } - \frac{\mu}{2\Gamma} \norm{x_f^k - x_{g}^{k}}_{\hat{D}_k}^{2} \nonumber\\&+ \left(\frac{L\gamma_k}{2e} - 1\right)\gamma_k\norm{\nabla f(x_g^k)}^2_{\hat{D}_{k}^{-1}}.
\end{align}
\begin{align}
\label{eq: n-third}
    f(x_f^{k + 1}) &\leq f(x^*) - \ev{ \nabla f(x_{g}^{k}), x^* - x_{g}^{k} } - \frac{\mu}{2\Gamma} \norm{x^* - x_{g}^{k}}_{\hat{D}_k}^{2} \nonumber\\&+ \left(\frac{L\gamma_k}{2e} - 1\right)\gamma_k\norm{\nabla f(x_g^k)}^2_{\hat{D}_{k}^{-1}}.
\end{align}
Summing \eqref{eq: n-sec} and \eqref{eq: n-third} with multiplicative factors $2\gamma_k \xi_k(\xi_k - 1)$ and $2\gamma_k \xi_k$ respectively, and substituting the result with \eqref{eq: n-first} into \eqref{eq: descent-nest-start}, we have
\begin{align*}
    \norm{x^{k+1} - x^{*}}_{\hat{D}_{k}}^{2} &+ 2\gamma_k\xi_k^2 f(x_f^{k+1}) \nonumber\\&\leq \left(1 - \frac{1}{\xi_k}\right)\norm{x^k - x^*}^2_{\hat{D}_k} + \left(\frac{1}{\xi^2} -\frac{1}{\xi}\right)\norm{x^k - x_g^k}^2_{\hat{D}_k} \nonumber\\&+ \frac{1}{\xi_k}\norm{x_g^k - x^*}^2_{\hat{D}_k} + \gamma_{k}^{2}\xi_{k}^{2}\norm{\nabla f(x_{g}^{k})}_{\hat{D}_{k}^{-1}}^{2} + 2\gamma_k \xi_k f(x^*) \nonumber\\&- \frac{\gamma_k\xi_k \mu}{\Gamma}\norm{x^* - x_{g}^{k}}_{\hat{D}_k}^{2} + 2\gamma_k \xi_k(\xi_k - 1) f(x_f^k) \nonumber\\&- \frac{\gamma_k \xi_k(\xi_k - 1)\mu}{\Gamma} \norm{x_f^k - x_{g}^{k}}_{\hat{D}_k}^{2} + 2\gamma_k^2\xi_k^2 \left(\frac{L\gamma_k}{2e} - 1\right)\norm{\nabla f(x_g^k)}^2_{\hat{D}_{k}^{-1}}.
\end{align*}
Combining the terms, we have
\begin{align*}
    \norm{x^{k+1} - x^{*}}_{\hat{D}_{k}}^{2} &+ 2\gamma_k\xi_k^2 f(x_f^{k+1}) \nonumber\\&\leq \left(1 - \frac{1}{\xi_k}\right)\norm{x^k - x^*}^2_{\hat{D}_k} + \left(\frac{1}{\xi^2} -\frac{1}{\xi}\right)\norm{x^k - x_g^k}^2_{\hat{D}_k} \nonumber\\&+ \left(\frac{1}{\xi_k} - \frac{\gamma_k\xi_k \mu}{\Gamma}\right)\norm{x_g^k - x^*}^2_{\hat{D}_k} + 2\gamma_k \xi_k f(x^*) \nonumber\\&+ 2\gamma_k \xi_k(\xi_k - 1) f(x_f^k) - \frac{\gamma_k \xi_k(\xi_k - 1)\mu}{\Gamma} \norm{x_f^k - x_{g}^{k}}_{\hat{D}_k}^{2} \nonumber\\&+ \gamma_k^2\xi_k^2 \left(\frac{L\gamma_k}{e} - 1\right)\norm{\nabla f(x_g^k)}^2_{\hat{D}_{k}^{-1}},
\end{align*}
where the scalar products are eliminated by multiplicative factors mentioned before. Subtracting $2\gamma_k\xi_k^2 f(x^*)$ from both sides, and using that $\xi_k \geq 1, \frac{\xi_k^2\gamma_k\mu}{\Gamma} \geq 1$ and $\gamma_k \leq \frac{e}{L}$, we obtain
\begin{align*}
    \norm{x^{k+1} - x^{*}}_{\hat{D}_{k}}^{2} &+ 2\gamma_k\xi_k^2 (f(x_f^{k+1}) - f(x^*)) \\&\leq \left(1 - \frac{1}{\xi_k}\right)\norm{x^k - x^*}^2_{\hat{D}_k} + 2\gamma_k\xi_k(\xi_k - 1) (f(x_f^{k}) - f(x^*))\\&= \left(1 - \frac{1}{\xi_k}\right)\left[\norm{x^{k} - x^{*}}_{\hat{D}_{k}}^{2} + 2\gamma_k\xi_k^2 (f(x_f^{k}) - f(x^*))\right],
\end{align*}
what finishes the proof.  
\end{proof} 
\noindent To obtain the convergence, let us formulate a final theorem.
\begin{tcolorbox}[colback=gray!25, colframe=gray!10, sharp corners]
\begin{theorem}
Suppose that Assumptions \ref{asm: strong}, \ref{asm: smooth} and \ref{asm: scaling-matr}  hold. Then, after $K$ iterations of \cref{alg:scaling_algorithm_nesterov} with $\gamma_k \equiv \frac{e}{L}, \xi_k \equiv \sqrt{\frac{L \Gamma}{\mu e}}$ and $\theta_k \equiv \frac{\sqrt{L\Gamma}}{\sqrt{\mu e} + \sqrt{L\Gamma}}$, we have
\begin{align*}
    \norm{x^{k+1} - x^{*}}_{\hat{D}_{k}}^{2} &\leq \exp\left(-K\sqrt{\frac{\mu e}{4L\Gamma}}\right)\left[\norm{x^{0} - x^{*}}_{\hat{D}_{-1}}^{2} + \frac{2\Gamma}{\mu} (f(x^{0}) - f(x^*))\right].
\end{align*}
\end{theorem}
\end{tcolorbox}
\begin{proof}
Starting with \cref{lem: descent-nesterov}:
\begin{align*}
    \norm{x^{k+1} - x^{*}}_{\hat{D}_{k}}^{2} &+ 2\gamma_k\xi_k^2 (f(x_f^{k+1}) - f(x^*)) \\&\leq \left(1 - \frac{1}{\xi_k}\right)\left[\norm{x^{k} - x^{*}}_{\hat{D}_{k}}^{2} + 2\gamma_k\xi_k^2 (f(x_f^{k}) - f(x^*))\right].
\end{align*}
We use the \cref{prop: ultra} with $\beta_2 \geq 1 - \frac{1}{2C}\sqrt{\frac{\mu e}{L\Gamma}}$ such that
\begin{align*}
    \norm{x^{k} - x^{*}}_{\hat{D}_{k}}^{2} \leq \left(1 + \frac{1}{2\xi_k}\right)\norm{x^{k} - x^{*}}_{\hat{D}_{k-1}}^{2}.
\end{align*}
Consequently, with $(1 - x)\left(1 + \frac{x}{2}\right) \leq \left(1  - \frac{x}{2}\right)$ we get
\begin{align*}
    \norm{x^{k+1} - x^{*}}_{\hat{D}_{k}}^{2} &+ 2\gamma_k\xi_k^2 (f(x_f^{k+1}) - f(x^*)) \\&\leq \left(1 - \frac{1}{2\xi_k}\right)\left[\norm{x^{k} - x^{*}}_{\hat{D}_{k-1}}^{2} + 2\gamma_k\xi_k^2 (f(x_f^{k}) - f(x^*))\right].
\end{align*}
After substituting $\gamma_k, \xi_k$ and $\theta_k$, with $(1 - x)\leq \exp(-x)$, the recursion provides a final bound.  
\end{proof}
\begin{remark}
    For Theorem 2 it is enough to apply \cref{prop: scale-bound} to obtain the convergence in $L_2$-space.
\end{remark}
\end{document}